\documentclass[12pt,a4paper]{article}

\usepackage{graphicx}
\usepackage{amsmath,amsfonts}
\usepackage{mathrsfs}
\usepackage[T1]{fontenc}
\usepackage[utf8]{inputenc}
\usepackage{authblk}
\usepackage{enumerate}
\usepackage{capt-of}
\usepackage{calc}
\usepackage{float}
\allowdisplaybreaks[1]

\newtheorem{theorem}{Theorem}
\newtheorem{lemma}{Lemma}
\newtheorem{proposition}{Proposition}
\newtheorem{corollary}{Corollary}

\let\scr\mathscr

\def\Pb{\mathbf{P}}

\def\Ex{\mathbf{E}}

\def\AA{\mathbb{A}}

\def\CC{\mathbb{C}}
\def\DD{\mathbb{D}}

\def\TT{\mathbb{T}}

\def\RR{\mathbb{R}}
\def\MM{\mathbb{M}}

\def\II{\mathbb{I}}
\def\UU{\mathbb{U}}

\def\1{\mbox{1\hspace{-.25em}I}}
\newcommand{\Liminf}{\mathop{\underline{\lim}}\limits}

\begin{document}
\title{Poisson Source Localization on the Plane. Smooth Case.}
\author[1]{O.V. Chernoyarov}
\author[2]{Yu.A. Kutoyants}

\affil[1,2]{\small  National Research University ``MPEI'', Moscow, Russia}
\affil[2]{\small  Le Mans University,  Le Mans,  France}
\affil[1,2]{\small  Tomsk State University, Tomsk, Russia}

\maketitle

\begin{abstract}
We consider the problem of localization of Poisson source by the
observations of inhomogeneous Poisson processes. We suppose that there are $k$
detectors on the plane and each detector provides the observations of Poisson
processes whose intensity functions depend on the position of the 
emitter. We describe the properties of the maximum likelihood and Bayesian
estimators. We show that under regularity conditions these estimators are
consistent, asymptotically normal and asymptotically efficient. Then we
propose some simple consistent estimators and this estimators are further used
to construct asymptotically efficient One-step MLE-process. 
\end{abstract}

\bigskip{} \textbf{Key words}: Inhomogeneous Poisson process,
source localization , GPS-localization, sensors, maximum likelihood estimator, Bayes
estimators,  one-step MLE-process.
 \bigskip{}
\date{}

\section{Introduction}

We consider the problem of estimation of the position $\vartheta
=\left(x_0,y_0\right)$ of a source emitting Poisson
signals which are  received  by distributed on the plane $k$ sensors \cite{Kn10}. We
suppose that the source starts emission at the instant $t=0$ and the $j$-th
sensor receives the data, which can be described as inhomogeneous Poisson
process $X_j=\left(X_j\left(t\right), 0\leq t\leq T\right)$, whose intensity
function $\lambda _j\left(\vartheta_0 ,t\right)=\lambda _{j}\left(t-\tau
_j\right)+\lambda _0, 0\leq t\leq T$ increases at the moment $t=\tau _j$ of
arriving the signal.  Here $\lambda _0>0$ is the intensity of the Poisson
noise and $\tau _j$ is the time needed for the signal to arrive at the $j$-th
detector. For the $j$-th detector localized at the point $\vartheta
_j=\left(x_j,y_j\right)$ we have $\tau _j\left(\vartheta _0\right)=\nu
^{-1}\left\|\vartheta _j-\vartheta _0\right\|$, where $\nu >0$ is the known
rate of propagation of the signal and $\left\|\cdot \right\|$ is Euclidean
norm on the plane. We suppose that $\lambda _j\left(t\right)=0$ for $t\leq 0$.
Therefore we have $k$ independent inhomogeneous Poisson processes
$X=\left(X_1,\ldots,X_k\right)$ with intensities depending on $\tau
_j\left(\vartheta _0\right)$. We suppose that the position of the source
$\vartheta _0\in\Theta $ is unknown and we have to estimate $\vartheta _0$ by
the observations $X=\left(X_1,\ldots,X_k\right)$. Here $\Theta \subset {\cal
  R}^2$ is a convex bounded set.

Note that the same mathematical model we have in the problem of
GPS-localization on the plane \cite{XL13}. Indeed, in this case we have $k$
emitters with known 
positions and an object which receives these signals and has to estimate its
own position. Therefore, we have observations of $k$ inhomogeneous Poisson
processes with the intensity functions depending on the position of the object
and we have to estimate the coordinates of this object. 

Due to importance of such type of models  in many applied problems there
exists a wide literature devoted to the different algorithms of localization
 (see the introduction in the work \cite{FKT18} and the references
there in). It seems that  the mathematical study of this class of models was
not yet sufficiently developed.  The statistical models of inhomogeneous
Poisson processes with intensity functions having discontinuities along some
curves depending on unknown parameters were considered in \cite{Kut98},
Sections 5.2 and 5.3. Statistical inference for point processes can be found
in the works \cite{Ka91}, \cite{SM91} and \cite{St10}.

 We are interested in the models of observations which
allow the estimation with the small errors: $\Ex_{\vartheta_0 }\left(\bar
\vartheta -\vartheta _0\right)^2=o\left(1\right) $. As usual in such situations
as we said ``small error'' we have to consider some asymptotic statement. 
 The small errors can be obtained, for example, 
if the intensity of the signal takes large values or we have periodical
Poisson process. Another possibility is to have many sensors.  We take the
model with large intensity functions $\lambda_j \left(\vartheta_0,
t\right)=\lambda_{j,n} \left(\vartheta_0, t\right) $, which can be written as
follows
\begin{align*}
\lambda_{j,n} \left(\vartheta_0, t\right)=n\lambda_j \left(t-\tau _j\right)+n\lambda
_0,\qquad 0\leq t\leq T
\end{align*}
or in equivalent form
\begin{align*}
\lambda_{j,n} \left(\vartheta_0, t\right)=n\lambda_j \left(t-\tau
_j\right)\1_{\left\{t\geq \tau _j\left(\vartheta _0\right)\right\}}+n\lambda 
_0,\qquad 0\leq t\leq T.
\end{align*}
Here $n$ is a ``large parameter'' and we study estimators as $n\rightarrow
\infty $. For example, such model we can be obtained if we have  $k$ clusters and
in each cluster we have $n$ detectors.

The likelihood ratio function $L\left(\vartheta ,X^n\right) $  is
\begin{align*}
\ln L\left(\vartheta ,X^n\right)&=\sum_{j=1}^{k}\int_{\tau _j}^{T}\ln \left(1+\frac{\lambda_j
  \left(t-\tau _j\right)}{\lambda _0}\right) {\rm
  d}X_j\left(t\right) -n\sum_{j=1}^{k}\int_{\tau _j}^{T}\lambda_j
  \left(t-\tau _j\right){\rm d}t. 
\end{align*}
Here $\tau _j=\tau _j\left(\vartheta \right)$ and $X=\left(X_j\left(t\right),0\leq
t\leq T, j=1, \ldots,k\right)$ are counting processes from  $k$ detectors. Having this
likelihood ratio formula we define the  maximum likelihood estimator (MLE)
$\hat\vartheta _n$   and Bayesian estimator (BE)  
$\tilde\vartheta _n$ by the ``usual'' relations
\begin{align}
\label{mle}
L\left(\hat\vartheta _n ,X^n\right)=\sup_{\vartheta \in \Theta }L\left(\vartheta ,X^n\right)
\end{align}
and
\begin{align}
\label{be}
\tilde\vartheta _n=\frac{\int_{\Theta }^{}\vartheta  p\left(\vartheta
  \right)L\left(\vartheta ,X^n\right){\rm d}\vartheta }{\int_{\Theta }^{}  p\left(\vartheta
  \right)L\left(\vartheta ,X^n\right){\rm d}\vartheta }.
\end{align}
Here $p\left(\vartheta \right),\vartheta \in \Theta $ is the prior
density. We suppose that it is  positive, continuous function on $\Theta $. If
the equation \eqref{mle} has more than one solution then any of these
solutions can be taken as the MLE.
 In
the section 3 we  consider another consistent  estimator.

There are several types  of  statistical  problems depending on the
regularity of the function $\lambda_j \left(\cdot  \right)$. In particularly, the
rate of convergence of the mean square error of the estimator $\bar\vartheta
_n$ is
\begin{align*}
\Ex_{\vartheta_0 }\left(\bar
\vartheta_n -\vartheta _0\right)^2=\frac{C}{n^\gamma }\left(1+o\left(1\right)\right),
\end{align*}
where the parameter $\gamma >0$ depends on the regularity of the function
$\lambda \left(\cdot \right)$.

Let us recall three of them using the following intensity functions
\begin{align}
\label{int}
\lambda _{j,n}\left(\vartheta _0,t\right)=an\left|t-\tau _j\left(\vartheta
_0\right)\right|^\kappa \1_{\left\{t\geq \tau _j\left(\vartheta
  _0\right)\right\}}+n\lambda _0, \quad 0\leq t\leq T.
\end{align}
We suppose that $ a>0, \lambda _0>0$ and known, the set $\Theta $ is such that
for all $\vartheta\in\Theta $ the instants $\tau _j\left(\vartheta \right)\in
\left(0,T\right)$.

\bigskip
{{\bf Here is Fig. 1}}

\bigskip

\begin{figure}[ht]
\hspace{1cm}\includegraphics[width=12cm]   {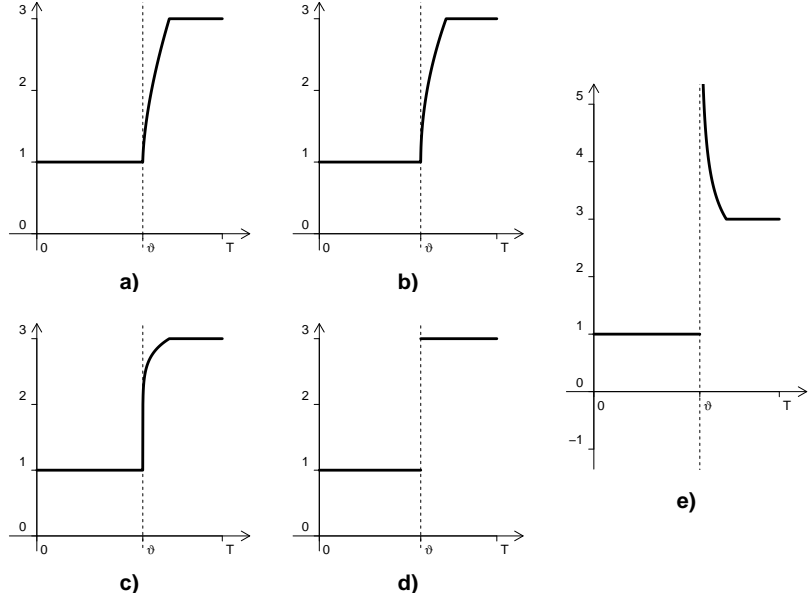}
\caption{Intensity functions: {\bf a}) $\kappa =\frac{5}{8}$, {\bf b}) $\kappa
  =\frac{1}{2}$, {\bf c}) $\kappa =\frac{1}{8}$, {\bf d}) $\kappa =0$, {\bf e})
  $ \kappa =-\frac{3}{8}$. }
\end{figure}

\begin{description}
\item[{\bf a}) Smooth case.] Suppose that the  $\kappa >\frac{1}{2} $, then the problem
  of parameter estimation is regular, the estimators are asymptotically normal
  and 
\begin{align*}
\Ex_{\vartheta_0 }\left\|\tilde
\vartheta_n -\vartheta _0\right\|^2=\frac{C}{n }\left(1+o\left(1\right)\right),\qquad \gamma =1.
\end{align*}

\item[{\bf b}) Smooth case.] If $\kappa =\frac{1}{2}$, then
\begin{align*}
\Ex_{\vartheta_0 }\left\|\tilde
\vartheta_n -\vartheta _0\right\|^2=\frac{C}{n\ln n
}\left(1+o\left(1\right)\right).
\end{align*}

  \item[{\bf c }) Cusp-type case.] This case is intermediate between the
    smooth and change-point cases. Suppose that $\kappa \in
    \left(0,\frac{1}{2}\right)$. Then
\begin{align*}
\Ex_{\vartheta_0 }\left\|\tilde
\vartheta_n -\vartheta _0\right\|^2=\frac{C}{n^{\frac{2}{2\kappa +1}}
}\left(1+o\left(1\right)\right),\qquad \gamma =\frac{2}{2\kappa +1}. 
\end{align*}

\item[{\bf  d}) Change point  case.] Suppose that $\kappa =0$. Then 
\begin{align*}
\Ex_{\vartheta_0 } \left\|\tilde \vartheta_n -\vartheta _0\right\|
^2=\frac{C}{n^2 }\left(1+o\left(1\right)\right),\qquad \gamma =2.
\end{align*}
\item[{\bf e }) Explosion case.]  Suppose that $\kappa \in
    \left(-\frac{1}{2},0\right)$. Then
\begin{align*}
\Ex_{\vartheta_0 }\left\|\tilde
\vartheta_n -\vartheta _0\right\|^2=\frac{C}{n^{\frac{2}{2\kappa +1}}
}\left(1+o\left(1\right)\right),\qquad \gamma =\frac{2}{2\kappa +1}. 
\end{align*}
\end{description}

The smooth case {\bf a}) is studied in this work. See as well the work
\cite{BMBD15}, where the similar model was considered.  The case {\bf b}) is
discussed below in the Section 4. For the Cusp case {\bf c}) see \cite{D03},
\cite{CDK18}. 
The change-point    case {\bf d}) is  studied in
 \cite{FKT18}. The Explosion case {\bf d}) can be done using the technique
 developed in \cite{D11}.

\section{Main result}

Suppose that there exists a  source at some point $\vartheta
_0=\left(x_0,y_0\right)\in \Theta \subset {\cal R}^2$ and $k\geq 3$ sensors
(detectors) on the same plane located at the points $\vartheta
_j=\left(x_j,y_j\right), j=1,\ldots,k$. The source was activated at the
(known) instant $t=0$ and the  signals from the source (inhomogeneous Poisson
processes) are registered by all  $k$  detectors.  The signal arrives at the
$j$-th detector at the instant $ \tau 
_j$. Of course, $\tau _j=\tau _j\left(\vartheta _0\right)$ is the  time
necessary for the signal to arrive in the $j$-th detector defined by the
relation
\begin{align*}
\tau _j\left(\vartheta _0\right)=\nu ^{-1}\left\|\vartheta _j-\vartheta
_0\right\|, 
\end{align*}
where $\nu >0$ is the known speed of propagation of the signal and
$\left\|\cdot \right\|$ is the Euclidean norm (distance) in ${\cal R}^2$.

The intensity function of the Poisson process $X_j^n=\left(X_j\left(t\right), 
0\leq t\leq T\right)$ registered by the $j$-th detector is
\begin{align*}
\lambda _j\left(\vartheta ,t\right)=n\lambda _j\left(t-\tau _j
\right)+n\lambda _0,\qquad 0\leq t\leq T .
\end{align*}
Here $n\lambda _j\left(t-\tau _j\right) $ is the intensity function of the
signal and $n\lambda _0>0 $ is the intensity of the noise. For simplicity of the
exposition we suppose that the noise level in all detectors is the same. 

Introduce the notations:
\begin{align*}
\alpha _j&=\inf_{\vartheta \in\Theta }\tau _j\left(\vartheta \right),\quad
\beta _j=\sup_{\vartheta \in\Theta }\tau _j\left(\vartheta \right),\qquad
j=1,\ldots,k\\ 
J_j\left(\vartheta \right)&=\frac{1}{\nu^2\left\|\vartheta
  _j-\vartheta \right\|^2} \int_{\tau _j\left(\vartheta \right)}^{T}\frac{
  \lambda _j'\left(t-\tau _j\left(\vartheta \right)\right)^2}{\lambda
  _j\left(t-\tau _j\left(\vartheta \right)\right)+\lambda _0}\; {\rm
  d}t,\\ 
\langle a,b\rangle_\vartheta &=\sum_{j=1}^{k}a_jb_jJ
_j\left(\vartheta \right),\qquad \left\|a\right\|_\vartheta
^2=\sum_{j=1}^{k}a_j^2J _j\left(\vartheta \right).
\end{align*}
Recall that $\lambda _j'\left(t-\tau _j\right)=0 $ for $0\leq t\leq \tau _j $
and note that $\langle a,b\rangle_\vartheta$ and $\left\|a\right\|_\vartheta $
are formally the scalar product and the norm in ${\cal R}^k$ of the vectors
$a=\left(a_1,\ldots,a_k\right)^\TT,$ $b=\left(b_1,\ldots,b_k\right)^\TT $ with
weights $\rho _j$ but the both depend on $\vartheta $ by a very special way.
The Fisher information matrix $\II_n\left(\vartheta
\right)=n\II\left(\vartheta \right) $, where $\vartheta =\left(x,y\right)$ and
\[\II\left(\vartheta \right)=\left(
\begin{array}{cc}
\left\|x-x _0\right\|_\vartheta ^2, &\langle
\left(x-x_0\right),\left(y-y_0\right)\rangle_\vartheta
\\
\langle \left(x-x_0\right),\left(y-y_0\right)\rangle_\vartheta,&\left\|y-y
_0\right\|_\vartheta ^2 \\  
\end{array}
\right).\]
Here $x=\left(x_1,\ldots,x_k\right), y=\left(y_1,\ldots,y_k\right)$ and
$x_0=\left(x_0,\ldots,x_0\right)$ {\it etc.}

Further, we suppose that $\beta _j<T$ and that the functions $\lambda
_j\left(t\right), j=1,\ldots,k$ are defined on the sets ${\cal T}_j=\left[-\beta
  _j,T-\alpha _j\right]$.

\bigskip

Regularity conditions  ${\scr R}$.\\
${\scr R}_1$. {\it  For all $j=1,\ldots,k$ the functions   }
$$
\lambda _j\left(t\right)=0,\quad  t\in 
\left[-\beta _j,0\right],\qquad   {and}\qquad \lambda
_j\left(t\right)> 0 ,\quad  t\in (0,T-\alpha _j] 
$$
${\scr R}_2$. {\it  The functions   $\lambda _j\left(t\right),t\in {\cal T}_j, 
j=1,\ldots,k$ have two  continuous derivatives $\lambda_j
'\left(\cdot \right)$ and  $\lambda_j ''\left(\cdot \right)$}.\\
${\scr R}_3$. {\it The Fisher information matrix is uniformly non degenerate}
\begin{align*}
\kappa _1=\inf_{\vartheta \in\Theta
}\inf_{\left|e\right|=1}e^\TT\II\left(\vartheta \right)e>0.
\end{align*}
${\scr R}_4$. {\it There are at least three detectors which are not on the same
  line.}

\bigskip

Remark, that if all detectors are on the same line, then the 
consistent identification is impossible because the same signals come from the
symmetric with respect to this line  possible locations of the source. 

According to Lemma \ref{L1} below the family of measures
$\left(\Pb_\vartheta^{\left(n\right)} ,\vartheta \in\Theta \right)$ induced by
the Poisson processes $X^n=\left(X_1^n,\ldots,X_k^n\right)$ in the space of
their realizations is {\it locally asymptotically normal} and therefore we
have the following minimax Hajek-Le Cam's lower boud on the mean square errors
of all estimators $\bar\vartheta _n$: for any $\vartheta _0\in\Theta $
\begin{align*}
\lim_{\delta \rightarrow 0}\Liminf_{n\rightarrow \infty
}\sup_{\left\|\vartheta -\vartheta _0\right\|\geq \delta } n\Ex_\vartheta
\left\|\bar\vartheta _n-\vartheta \right\|^2 \geq \Ex_{\vartheta _0}
\left\|\zeta \right\|^2 ,\qquad \zeta \sim {\cal
  N}\left(0,\II\left(\vartheta_0 \right)^{-1} \right).
\end{align*}
We call the estimator $\bar\vartheta _n$ asymptotically efficient, if for all
$\vartheta _0\in\Theta $ we have the equality
\begin{align*}
\lim_{\delta \rightarrow 0}\lim_{n\rightarrow \infty
}\sup_{\left\|\vartheta -\vartheta _0\right\|\geq \delta } n\Ex_\vartheta
\left\|\bar\vartheta _n-\vartheta \right\|^2 = \Ex_{\vartheta _0}
\left\|\zeta \right\|^2 .
\end{align*}
For the proof of this bound see, e.g., \cite{IH81}, Theorem 2.12.1.

\begin{theorem}
\label{T1} Let the conditions ${\cal R}$ be fulfilled then the MLE
$\hat\vartheta _n$ and BE $\tilde\vartheta _n$ are uniformly consistent,
asymptotically normal 
\begin{align*}
\sqrt{n}\left( \hat\vartheta _n-\vartheta _0 \right)\Longrightarrow {\cal
  N}\left(0,\II\left(\vartheta_0 \right)^{-1} \right),\qquad \sqrt{n}\left(
\tilde \vartheta _n-\vartheta _0 \right)\Longrightarrow {\cal
  N}\left(0,\II\left(\vartheta_0 \right)^{-1} \right),
\end{align*}
for any $p>0$
\begin{align*}
\lim_{n\rightarrow \infty } n^{\frac{p}{2}}\Ex_{\vartheta
  _0}\left\|\hat\vartheta _n-\vartheta _0 \right\|^p =\Ex_{\vartheta
  _0}\left\| \zeta \right\|^p,\qquad \lim_{n\rightarrow \infty }
n^{\frac{p}{2}}\Ex_{\vartheta _0}\left\|\tilde\vartheta _n-\vartheta _0
\right\|^p =\Ex_{\vartheta _0}\left\| \zeta \right\|^p,
\end{align*}
where $\zeta \sim {\cal N}\left(0,\II\left(\vartheta_0 \right)^{-1} \right) $
and the both estimators are asymptotically efficient.

\end{theorem}
{\bf Proof.} The proof of this theorem is based on two general results by
Ibragimov and Khasminskii \cite{IH81} presented  in the Theorems 1.10.1 and
1.10.2. We have to check the conditions of these theorems given in terms of
normalized likelihood ratio
\begin{align*}
Z_n\left(u\right)=\frac{L\left(\vartheta_0+\frac{u}{\sqrt{n}},X^n\right)
}{L\left(\vartheta_0 ,X^n\right)},\qquad u\in
\UU_n=\left\{u:\;\vartheta_0+\frac{u}{\sqrt{n}}\in \Theta  \right\}.
\end{align*}
Introduce the limit likelihood ratio
\begin{align*}
Z\left(u\right)=\exp\left\{\langle u,\Delta\left(\vartheta _0\right) \rangle
-\frac{1}{2}u^\TT\II\left(\vartheta_0 \right)u\right\} ,\qquad u\in {\cal
  R}^2 .
\end{align*}
Here $\Delta\left(\vartheta _0\right)\sim {\cal N}\left(0,\II\left(\vartheta_0
\right)\right) $. 

Suppose that we already proved the weak  convergence 
\begin{align*}
Z_n\left(\cdot \right)\Longrightarrow Z\left(\cdot \right).
\end{align*}
Then the limit distributions  of the mentioned estimators are obtained as
follows (see \cite{IH81}). Below we change the variables $\vartheta =\vartheta
_0+\frac{u}{\sqrt{n}}$  and ${\cal B}\subset {\cal R}^2$ is a bounded set.

For the MLE we have
\begin{align*}
&\Pb_{\vartheta _0}\left(\sqrt{n}{\left(\hat\vartheta _n-\vartheta_0\right) }
  \in {\cal B}\right)\\ 
&=\Pb_{\vartheta _0}\left\{ \sup_{\sqrt{n}\left(\vartheta
    -\vartheta_0\right)\in{\cal B}} L\left(\vartheta
  ,X^T\right)>\sup_{\sqrt{n}\left(\vartheta -\vartheta_0\right)\in{\cal B}^c}
  L\left(\vartheta ,X^T\right) \right\} \\ 
&=\Pb_{\vartheta _0}\left\{
  \sup_{\sqrt{n}\left(\vartheta 
    -\vartheta_0\right)\in{\cal B}} \frac{L\left(\vartheta
  ,X^T\right)}{L\left(\vartheta_0
  ,X^T\right) }>\sup_{\sqrt{n}\left(\vartheta -\vartheta_0\right)\in{\cal B}^c}
  \frac{L\left(\vartheta ,X^T\right)}{L\left(\vartheta_0 ,X^T\right) }
  \right\}\\
 & =\Pb_{\vartheta _0}\left\{ \sup_{u\in {\cal B}
    ,u\in\UU_n } Z_n \left(u\right)>\sup_{u \in{\cal B}^c,u\in\UU_n } Z_n
  \left(u\right) \right\} \\ &\qquad \qquad \longrightarrow
  \Pb_{\vartheta _0}\left\{ \sup_{u \in{\cal B}} Z\left(u\right)>\sup_{u\in{\cal B}^c }
  Z\left(u\right) \right\} =\Pb_{\vartheta _0}\left(\zeta \in {\cal B}\right)  .
\end{align*}
It is easy to see that $\zeta =\arg\max_u Z\left(u\right)$

For the BE we have (once more we change the
variables $\theta_u=\vartheta _0+\frac{u}{\sqrt{n}} $): 
\begin{align*}
\tilde\vartheta _n &=\frac{\int_{\Theta }^{}\theta p\left(\theta
  \right)L\left(\theta ,X^T\right){\rm d}\theta }{\int_{\Theta }^{}
  p\left(\theta \right)L\left(\theta ,X^T\right){\rm d}\theta}=\vartheta
_0+\frac{1}{\sqrt{n}} \frac{\int_{\UU_n }^{}u p\left(\theta_u
  \right)L\left(\theta_u ,X^T\right){\rm d}u }{\int_{\UU_n }^{}
  p\left(\theta_u \right)L\left(\theta_u ,X^T\right){\rm d}u}
\\ &=\vartheta
_0+\frac{1}{\sqrt{n}} \frac{\int_{\UU_n }^{}u p\left(\theta_u \right)Z_n
  \left(u\right){\rm d}u }{\int_{\UU_n }^{} p\left(\theta_u \right)Z_n
  \left(u\right){\rm d}u}.
\end{align*}
Hence
\begin{align*}
\sqrt{n}\left(\tilde\vartheta _n-\vartheta _0\right)=\frac{\int_{\UU_n }^{}u
  p\left(\theta_u \right)Z_n \left(u\right){\rm d}u }{\int_{\UU_n }^{}
  p\left(\theta_u \right)Z_n \left(u\right){\rm d}u}\Longrightarrow
\frac{\int_{{\cal R}^2 }^{}u Z\left(u\right){\rm d}u }{\int_{{\cal R}^2 }^{}
  Z\left(u\right){\rm d}u}=\zeta .
\end{align*}
Recall that $p\left(\theta _u\right)\rightarrow p\left(\vartheta
_0\right)>0$ and note that
\begin{align*}
\int_{{\cal R}^2
  }^{}u Z\left(u\right){\rm d}u& =\zeta\int_{{\cal R}^2
  }^{} Z\left(u\right){\rm d}u.
 \end{align*}
The properties of the $Z_n\left(u\right)$ required in the Theorems 1.10.1 and
1.10.2 \cite{IH81} are checked in the three lemmas below.  Remind that this
approach to the study of the properties of these estimators was applied  in
\cite{Kut79},  \cite{Kut98}. Here we use some obtained there inequalities.  

Introduce the vector of partial derivatives
\begin{align*}
\Delta _n\left(\vartheta _0,X^n\right)=\frac{1}{\sqrt{n}}\left( \frac{\partial
  \ln L\left(\vartheta_0 ,X^n\right)}{\partial 
  x_0},\, \frac{\partial  \ln L\left(\vartheta_0 ,X^n\right)}{\partial
  y_0} \right)^\TT .
\end{align*}

The convergence of finite dimensional distributions of the random field
$Z_n\left(u\right), u\in\UU_n$ to the finite dimensional distributions of the
limit random field $Z\left(u\right),u\in{\cal R}^2$ follows from the Lemma
\ref{L1} below.

\begin{lemma}
\label{L1} Let the conditions ${\scr R}$ be fulfilled,  then the family of
measures $\left\{ \Pb_\vartheta ,\vartheta \in\Theta \right\}$ is locally
asymptotically normal (LAN), i.e., the random process $Z_n\left(u\right),u\in
\UU_n$ for any $\vartheta _0\in\Theta $ admits the representation
\begin{align}
\label{lan1}
 Z_n\left(u\right)=\exp\left\{\langle u,\Delta _n\left(\vartheta
_0,X^n\right)\rangle -\frac{1}{2}u^\TT\II\left(\vartheta_0 \right)u+r_n\right\},
 \qquad u\in \UU_n,
\end{align}
where the vector
\begin{align}
\label{lan2}
\Delta _n\left(\vartheta _0,X^n\right)\Longrightarrow \Delta \left(\vartheta
_0\right) \sim {\cal
  N}\left(0,\II\left(\vartheta_0 \right) \right) 
\end{align}
and $r_n\rightarrow 0$. 
\end{lemma}
{\bf Proof.}  Let us denote $\lambda _j\left(t,u\right)=\lambda _j\left(t-\tau
_j (\theta _u \right))$ and put ${\rm d}\pi _{j,n}\left(t\right)= {\rm
  d}X_j\left(t\right) -n\left[\lambda _j\left(t-\tau _j\left(\vartheta
  _0\right)\right) +\lambda _0\right]{\rm d}t$. Then we can write
\begin{align*}
&\ln Z_n\left(u\right)=\sum_{j=1}^{k}\int_{0}^{T}\ln\left(\frac{\lambda
  _j\left(t,u)+\lambda _0\right)
}{\lambda _j\left(t,0\right)+\lambda _0
}\right){\rm d}\pi _{j,n}\left(t\right)\\
&\quad -n\sum_{j=1}^{k}\int_{0}^{T}\left[\frac{\lambda
  _j\left(t,u\right)+\lambda _0}{ \lambda _j\left(t,0\right)+\lambda _0 } 
  -1- \ln\left(\frac{\lambda
  _j(t,u)+\lambda _0
}{\lambda _j\left(t,0\right)+\lambda _0
}\right) \right] \left[\lambda _j\left(t,0\right)+\lambda _0 \right]{\rm d}t.
\end{align*}
Using the Taylor formula we obtain the relations
\begin{align*}
&\tau _j\left(\vartheta _0+\frac{u}{\sqrt{n}}\right)=\tau _j\left(\vartheta
  _0\right)- \frac{1}{\nu\sqrt{n} }\langle m_j,
  u\rangle+O\left(\frac{1}{n}\right), \\
 & m_j=\left(\frac{x_j-x_0}{\rho _j},\frac{y_j-y_0}{\rho _j} \right) ,\qquad
  \left\|m_j\right\|=1, \\
&\lambda _j(t-\tau
  _j(\vartheta _0+n^{-1/2}u))-\lambda _j(t-\tau
  _j(\vartheta _0))\\
&\qquad \quad \qquad =-n^{-1/2}{\lambda _j'(t-\tau
  _j(\vartheta _0))} \langle u, \frac{\partial \tau \left(\vartheta
  _0\right)}{\partial  \vartheta }\rangle+n^{-1}
  O\left({\left\|u\right\|^2}\right)\\
&\qquad \quad \qquad =n^{-1/2}\nu ^{-1}{\lambda _j'(t-\tau
  _j(\vartheta _0))} \langle m_j,  u\rangle+n^{-1}
  O\left({\left\|u\right\|^2}\right),\\
&\ln\left(\frac{\lambda  _j\left(t,u)+\lambda _0\right)
}{\lambda _j\left(t,0\right)+\lambda _0}\right)=\frac{\lambda _j'(t-\tau
  _j(\vartheta _0))}{\sqrt{n}\left[\lambda _j(t-\tau
  _j(\vartheta _0))+\lambda _0 \right]} \left(u, \frac{\partial \tau \left(\vartheta
  _0\right)}{\partial  \vartheta }\right)+
  O\left(\frac{1}{n}\right),\\
&\qquad \quad \qquad=\frac{\lambda _j'(t-\tau
  _j(\vartheta _0))}{\nu \sqrt{n}\left[\lambda _j(t-\tau
  _j(\vartheta _0))+\lambda _0 \right]} \langle m_j,  u\rangle+
  O\left(\frac{1}{n}\right),\\
&\frac{\lambda
  _j\left(t,u\right)+\lambda _0}{ \lambda _j\left(t,0\right)+\lambda _0 } 
  -1- \ln\left(\frac{\lambda
  _j(t,u)+\lambda _0
}{\lambda _j\left(t,0\right)+\lambda _0
}\right)\\
&\qquad \quad \qquad=\frac{1}{2n}\frac{\lambda _j'(t-\tau
  _j(\vartheta _0))^2}{\left[\lambda _j(t-\tau
  _j(\vartheta _0))+\lambda _0\right]^2 } \left(u, \frac{\partial \tau_j
    \left(\vartheta 
  _0\right)}{\partial  \vartheta }\right)^2+O\left(\frac{1}{n^{3/2}}\right)\\
&\qquad \quad \qquad=\frac{1}{2n\nu ^2}\frac{\lambda _j'(t-\tau
  _j(\vartheta _0))^2}{\left[\lambda _j(t-\tau
  _j(\vartheta _0))+\lambda _0\right]^2 }\langle m_j,
  u\rangle^2+O\left(\frac{1}{n^{3/2}}\right). 
\end{align*}
  Note that
\begin{align*}
\frac{\partial \tau _j\left(\vartheta _0\right)}{\partial
  x_0}=-\frac{x_j-x_0}{\nu \,\left\|\vartheta _j-\vartheta _0\right\|} ,\qquad
\frac{\partial \tau _j\left(\vartheta _0\right)}{\partial 
  y_0}=-\frac{y_j-y_0}{\nu \,\left\|\vartheta_j-\vartheta _0\right\|} 
\end{align*}
Therefore we can write
\begin{align*}
\frac{\partial \ln L\left(\vartheta_0 ,X^n\right)}{\partial
  x_0}=\sum_{j=1}^{k}\frac{\left(x_j-x_0\right)}{\nu
  \left\|\vartheta_j-\vartheta _0\right\|} \int_{\tau _j\left(\vartheta
  _0\right)}^{T}\frac{ \lambda 
  _j'\left(t-\tau _j\left(\vartheta _0\right)\right)}{\lambda _j\left(t-\tau
  _j\left(\vartheta _0\right)\right)+\lambda _0}\; {\rm d}\pi
_{j,n}\left(t\right). 
\end{align*}
 Hence
\begin{align*}
\Ex_{\vartheta _0}\left[\frac{\partial \ln L\left(\vartheta_0
    ,X^n\right)}{\partial
    x_0}\right]^2=n\sum_{j=1}^{k}\frac{\left(x_j-x_0\right)^2}
{\nu^2\left\|\vartheta_j-\vartheta 
  _0\right\|^2} \int_{\tau _j\left(\vartheta
  _0\right)}^{T}\frac{ \lambda _j'\left(t-\tau _j\left(\vartheta
  _0\right)\right)^2}{\lambda _j\left(t-\tau _j\left(\vartheta
  _0\right)\right)+\lambda _0}\; {\rm d}t
\end{align*}
and
\begin{align*}
\Ex_{\vartheta _0}\left[\frac{\partial \ln L\left(\vartheta_0 ,X^n\right)}{\partial
  x_0}\frac{\partial \ln L\left(\vartheta_0 ,X^n\right)}{\partial
  y_0}\right]=n\sum_{j=1}^{k}{\left(x_j-x_0\right)\left(y_j-y_0\right)}{}
J_j\left(\vartheta _0\right).
\end{align*}
These equalities justify the introduced above form of the Fisher information
matrix $\II\left(\vartheta _0\right)$. 

We have the representations
\begin{align*}
&\sum_{j=1}^{k}\int_{0}^{T}\ln\left(\frac{\lambda _j\left(t,u)+\lambda
    _0\right) }{\lambda _j\left(t,0\right)+\lambda _0 }\right){\rm d}\pi
  _{j,n}\left(t\right)\\ &\qquad =\frac{1}{\nu \sqrt{n}}\sum_{j=1}^{k} \langle
  m_j,u\rangle \int_{ \tau _j\left(\vartheta _0\right) }^{T}\frac{ \lambda
    _j'\left(t-\tau _j\left(\vartheta _0\right)\right)}{\lambda _j\left(t-\tau
    _j\left(\vartheta _0\right)\right)+\lambda _0}\; {\rm d}\pi
  _{j,n}\left(t\right)+o\left(1\right), \\ 
&\qquad= \langle u,\Delta
  _n\left(\vartheta _0,X^n\right) \rangle+o\left(1\right),
\end{align*}
and
\begin{align*}
 & n\sum_{j=1}^{k}\int_{0}^{T}\left[\frac{\lambda _j\left(t,u\right)+\lambda
      _0}{ \lambda _j\left(t,0\right)+\lambda _0 } -1- \ln\left(\frac{\lambda
      _j(t,u)+\lambda _0 }{\lambda _j\left(t,0\right)+\lambda _0 }\right )
    \right] \left[\lambda _j\left(t,0\right)+\lambda _0 \right]{\rm d}t\\
&\qquad =\frac{1}{2\nu ^2}\sum_{j=1}^{k} \langle m_j,u\rangle^2
  \int_{\tau _j\left(\vartheta
    _0\right)}^{T}\frac{ \lambda _j'\left(t-\tau _j\left(\vartheta
    _0\right)\right)^2}{\lambda _j\left(t-\tau _j\left(\vartheta
    _0\right)\right)+\lambda _0}\; {\rm d}t+o\left(1\right)\\
&\qquad =\frac{1}{2}\; u^\TT\,\II\left(\vartheta _0\right)u+o\left(1\right).
\end{align*}
Therefore we obtained \eqref{lan1}. To verify the convergence \eqref{lan2} we
introduce the vector $I_n=\left(I_{1,n},I_{2,n} \right) $, where
\begin{align*}
I_{1,n}&=\frac{1}{\sqrt{n}}\sum_{j=1}^{k}a_j\int_{\tau _j\left(\vartheta
    _0\right)}^{T}\frac{ \lambda _j'\left(t-\tau _j\left(\vartheta
    _0\right)\right)}{\lambda _j\left(t-\tau _j\left(\vartheta
    _0\right)\right)+\lambda _0}\; {\rm d}\pi _{j,n}\left(t\right),\\
I_{2,n}&=\frac{1}{\sqrt{n}}\sum_{j=1}^{k}b_j\int_{\tau _j\left(\vartheta
    _0\right)}^{T}\frac{ \lambda _j'\left(t-\tau
  _j\left(\vartheta 
    _0\right)\right)}{\lambda _j\left(t-\tau _j\left(\vartheta
    _0\right)\right)+\lambda _0}\; {\rm d}\pi _{j,n}\left(t\right),
\end{align*}
where the vectors $a,b\in{\cal R}^k$. Then the asymptotic normality of $I_n$
follows from the central limit theorem for stochastic integrals. See, e.g.,
Theorem 1.1 in \cite{Kut98}. Moreover, we have
\begin{align*}
I_{1,n}\Longrightarrow \sum_{j=1}^{k}a_j\int_{\tau _j\left(\vartheta
    _0\right)}^{T}\frac{ \lambda _j'\left(t-\tau _j\left(\vartheta
    _0\right)\right)}{{\lambda _j\left(t-\tau _j\left(\vartheta
    _0\right)\right)+\lambda _0}}\; {\rm d}W _j\left(\Lambda \left(\vartheta
_0,t\right)  \right), 
\end{align*}
where $W_j\left(\cdot \right), j=1,\ldots,k$ are independent Wiener processes
and 

 The conditions of this theorem can be easily
verified for the corresponding vectors $a,b$ given by the presentation $\Delta
_n\left(\vartheta ,X^n\right)$.

\begin{lemma} 
\label{L2}
 Let the condition ${\scr R}_2 $ be fulfilled,  then there exists a
  constant $C>0$, which does not depend on $n$ such that for any $R>0$
\begin{align}
\label{sm}
\sup_{\vartheta _0\in\Theta }\sup_{\left\|u_1\right\|+\left\|u_2\right\|\leq
  R}\left\|u_1-u_2\right\|^{-4} \Ex_{\vartheta _0} \left|
Z_n^{\frac{1}{4}}\left(u_1\right)-Z_n^{\frac{1}{4}}\left(u_2\right)\right|^4\leq 
C \left(1+R^2\right).
\end{align}
\end{lemma}
{\bf Proof.} The proof of this lemma follows from the proof of the Lemma 2.2
in \cite{Kut98} if we put there $m=2$. The difference between the models of
observations there and here is not essential for the presented there proof.

\begin{lemma} 
\label{L3}
Let the conditions ${\scr R} $ be fulfilled,  then there exists a
  constant $\kappa >0$, which does not depend on $n$ such that
\begin{align}
\label{pi}
\sup_{\vartheta _0\in\Theta }\Ex_{\vartheta _0} 
Z_n^{\frac{1}{2}}\left(u\right)\leq e^{-\kappa \left\|u\right\|^2}.
\end{align}
\end{lemma}
{\bf Proof.} Let us denote $\theta _u=\vartheta _0+\frac{u}{\sqrt{n}}$ and put 
\begin{align*}
Z_{j,n}\left(u\right)&=\exp\left\{\int_{0}^{T}\ln\left(\frac{\lambda
  _{j,n}\left(\theta_u,t\right)}{\lambda
  _{j,n}\left(\vartheta_0,t\right)}\right){\rm
  d}X_j\left(t\right)\right.\\
&\qquad \qquad \qquad \left.-\int_{0}^{T}\left[\lambda
  _{j,n}\left(\theta_u,t\right)-\lambda
  _{j,n}\left(\vartheta_0,t\right) \right]{\rm d}t\right\} .
\end{align*} 
Remind that (see Lemma 2.2  in \cite{Kut98})
\begin{align*}
\Ex_{\vartheta _0} Z_{j,n}^{\frac{1}{2}}\left(u\right)
=\exp\left\{-\frac{1}{2}\int_{0}^{T}\left[\sqrt{\lambda
  _{j,n}\left(\theta_u,t\right)}- \sqrt{\lambda
  _{j,n}\left(\vartheta_0,t\right)} \right]^2{\rm d}t       \right\} 
\end{align*}
Therefore we have the equality 
\begin{align}
\Ex_{\vartheta _0} Z_n^{\frac{1}{2}}\left(u\right)&=\prod_{j=1}^k
\Ex_{\vartheta _0} Z_{j,n}^{\frac{1}{2}}\left(u\right)\nonumber\\ 
&=\exp\left\{-\frac{1}{2}\sum_{j=1}^{k}\int_{0}^{T}\left[\sqrt{\lambda
  _{j,n}\left(\theta_u,t\right)}- \sqrt{\lambda
  _{j,n}\left(\vartheta_0,t\right)} \right]^2{\rm d}t       \right\} .
\label{32}
\end{align}
By Taylor formula for  $\left\|h\right\|\leq \delta $  we can write
\begin{align*}
\sum_{j=1}^{k}\int_{0}^{T}\left[\sqrt{\lambda
  _{j,n}\left(\vartheta_0+h,t\right)}- \sqrt{\lambda
  _{j,n}\left(\vartheta_0,t\right)} \right]^2{\rm d}t  =\frac{n}{4}\;
h^\TT\,\II\left(\vartheta _0\right)h \left(1+O\left(\delta \right)\right).
\end{align*}
Hence we can take such (small) $\delta>0 $ that for
$\frac{\left\|u\right\|}{\sqrt{n}}\leq \delta $  we have
\begin{align}
&\sum_{j=1}^{k}\int_{0}^{T}\left[\sqrt{\lambda
  _{j,n}(\vartheta_0+\frac{u}{\sqrt{n}},t)}- \sqrt{\lambda
  _{j,n}\left(\vartheta_0,t\right)} \right]^2{\rm d}t \geq \frac{1}{8}\;
u^\TT\,\II\left(\vartheta _0\right)u\nonumber\\
&\qquad \qquad =\frac{1}{8}\;
\frac{u^\TT}{\left\|u\right\|}\,\II\left(\vartheta
_0\right)\frac{u}{\left\|u\right\|}\;\left\|u\right\|^2\geq \frac{\kappa
  _1}{8}\;\left\|u\right\|^2, 
\label{33}
\end{align}
where $\kappa _1>0$ from the condition ${\scr R}_3$.

Let us denote
\begin{align*}
g\left(\delta \right)=\frac{1}{n}\inf_{\vartheta _0\in\Theta
}\inf_{\left\|\vartheta -\vartheta _0 
  \right\|>\delta }\sum_{j=1}^{k}\int_{0}^{T}\left[\sqrt{\lambda 
  _{j,n}\left(\vartheta,t\right)}- \sqrt{\lambda
  _{j,n}\left(\vartheta_0,t\right)} \right]^2{\rm d}t ,
\end{align*}
 and show that $g\left(\delta \right)>0$. Remark that $g\left(\delta \right)$
 does not depend on $n$. Indeed,
\begin{align*}
&\frac{1}{n}\int_{0}^{T}\left[\sqrt{\lambda 
  _{j,n}\left(\vartheta,t\right)}- \sqrt{\lambda
  _{j,n}\left(\vartheta_0,t\right)} \right]^2{\rm d}t\\
&\qquad =\frac{1}{n}\int_{0}^{T}\frac{\left[\lambda 
  _{j,n}\left(\vartheta,t\right)- \lambda
  _{j,n}\left(\vartheta_0,t\right)\right]^2}{\left[\sqrt{\lambda 
  _{j,n}\left(\vartheta,t\right)}+ \sqrt{\lambda
  _{j,n}\left(\vartheta_0,t\right)} \right]^2}{\rm d}t\\
&\qquad \geq \frac{1}{4\left(\lambda _M+\lambda
    _0\right)}\int_{0}^{T}{\left[\lambda  
  _{j}\left(t-\tau _j\left(\vartheta \right)\right)- \lambda
  _{j}\left(t-\tau _j\left(\vartheta_0 \right) \right)\right]^2}{\rm d}t.
\end{align*}
  Therefore if we suppose that $g\left(\delta \right)=0$ then this implies
  that there exists at least one point $\vartheta^*\in \Theta $ such that
  $\left\|\vartheta ^*-\vartheta _0\right\|\geq \delta $ and for all
  $j=1,\ldots,k$ we have
\begin{align*}
\int_{0}^{T}{\left[\lambda 
  _{j}\left(t-\tau _j\left(\vartheta^* \right)\right)- \lambda
  _{j}\left(t-\tau _j\left(\vartheta_0 \right) \right)\right]^2}{\rm d}t=0.
\end{align*}

Note that by the condition ${\scr R}_1$ we have the consistent estimation of
all ``delays'' $\tau _j.$ Indeed, the identifiability condition
\begin{align}
\label{idt}
g_j\left(\delta \right)\equiv \inf_{\left|\tau -\tau _0\right|\geq \delta
}\int_{0}^{T}{\left[\lambda  
  _{j}\left(t-\tau \right)- \lambda
  _{j}\left(t-\tau _0  \right)\right]^2}{\rm d}t>0
\end{align}
is fulfilled for all $j=1,\ldots,k$ and $\delta >0$. If for some $j$ and
$\delta >0$ we have $g_j\left(\delta \right)=0$, then this implies that there
exists $\tau _*\not=\tau _0$ such that the equality $\lambda
_{j}\left(t-\tau_* \right)= \lambda _{j}\left(t-\tau _0 \right) $ holds for
all $ t\in\left[0,T\right]$. This equality is impossible by the following
reason. Suppose that $\tau _0<\tau_* $, then on the interval $\left(\tau
_0,\tau _*\right)$ the function $\lambda _{j}\left(t-\tau _0 \right)>0 $ and
the function $\lambda _{j}\left(t-\tau _* \right)=0 $.

Therefore the condition \eqref{idt} is always fulfilled and we have the
consistency of the MLE $\hat\tau _{j,n}$ and BE $\tilde\tau _{j,n}$ of the
parameters $ \tau _j\left(\vartheta _0\right)$
\begin{align*}
\hat\tau _{j,n}\longrightarrow \tau _j\left(\vartheta
_0\right),\qquad \tilde\tau _{j,n}\longrightarrow \tau _j\left(\vartheta
_0\right),\qquad j=1,\ldots,k 
\end{align*}
and asymptotic normality of these estimators. For the proofs see
\cite{Kut98}. 

If $g\left(\delta \right)=0$, then there exist two points $\vartheta ^*$ and
$\vartheta _0$ such that we have  two sets $\left(\tau _j\left(\vartheta ^*\right),
j=1,\ldots,k\right)$ and $\left(\tau _j\left(\vartheta _0\right),
j=1,\ldots,k\right)$ with coinciding values $\tau _j\left(\vartheta
^*\right)=\tau _j\left(\vartheta _0\right) $ for all $j=1,\ldots,k$, i.e., the
distances
$\left\|\vartheta _j-\vartheta ^*\right\|$ and $\left\|\vartheta
_j-\vartheta_0\right\|$ for all $j=1,\ldots,k$ coincide. Note that such
situation is impossible due to geometric properties of the set of points
$\vartheta _1,\ldots,\vartheta _k$ satisfying the condition ${\scr R}_4$.

Hence, $g\left(\delta \right)>0$ and we can write for
$\frac{\left\|u\right\|}{\sqrt{n}} >\delta $
\begin{align}
\label{k}
\sum_{j=1}^{k}\int_{0}^{T}\left[\sqrt{\lambda
  _{j,n}(\vartheta_0+\frac{u}{\sqrt{n}},t)}- \sqrt{\lambda
  _{j,n}\left(\vartheta_0,t\right)} \right]^2{\rm d}t \geq ng\left(\delta
\right)\geq  g\left(\delta
\right)\frac{\left\|u\right\|^2}{ D^2},
\end{align}
where $D=\sup_{\vartheta _1,\vartheta _2\in\Theta
}\left\|\vartheta_1-\vartheta _2\right\|$.

 Let us denote $\kappa
_2=D^{-2}g\left(\delta \right)$ and put $\kappa =\min\left(\frac{\kappa
  _1}{16},\frac{\kappa _2}{2}\right)$, then from \eqref{32}, \eqref{33} and
\eqref{k} follows the estimate \eqref{pi}.   

\bigskip

The properties of the likelihood ratio field $Z_n\left(\cdot \right)$
established in the lemmas \ref{L1}-\ref{L2} are sufficient conditions for the
Theorems 1.10.1 and 1.10.2 in \cite{IH81}. Therefore the MLE $\hat\vartheta
_n$ and the BE $\tilde\vartheta _n$ have all mentioned in the Theorem \ref{T1}
properties. 

\section{Simple consistent estimator}

To find the MLE as solution of the equation \eqref{mle} can be computationally
quite complicate. It can be interesting to find some other estimators which
can be much more easily calculated. If we have just three detectors, then one
simple estimator was proposed by Pu \cite{Pu09}. 
 We consider the construction of the estimator in two steps.  First we
 solve $k$ one dimensional estimation problems of estimation of arriving times
 $\tau _j,j=1,\ldots,k$ and then having $k$ estimators $\hat\tau
 _{1,n},\ldots,\hat\tau _{k,n}$ we estimate the parameter $\vartheta _0$.

 We have $k$
inhomogeneous Poisson processes $X^n=\left(X_1^n,\ldots,X_k^n\right)$, where
$X_j^n=\left(X_j\left(t\right),0\leq t\leq T\right)$ is a Poisson process with
the intensity function
\begin{align*}
\lambda _{j,n}\left(\tau_j ,t\right)=n\lambda _j\left(t-\tau
_j\right)+n\lambda _0,\qquad 0\leq t\leq T .
\end{align*}
Therefore we have $k$ likelihood ratios
\begin{align*}
L\left(\tau_j ,X_j^n \right)=\exp\left\{\int_{\tau_j}^{T}\ln\left(1+\frac{\lambda
  _j\left(t-\tau_j \right)}{\lambda _0}\right){\rm
  d}X_j\left(t\right)-n\int_{\tau_j}^{T}\lambda _j\left(t-\tau _j\right){\rm
  d}t\right\}
\end{align*}
and can introduce $k$ MLEs $\hat\tau _{j,n},$ by the equations
\begin{align*}
L\left(\hat\tau _{j,n} ,X_j^n \right)=\sup_{\tau_j \in \Theta _j}L\left(\tau_j
,X_j^n \right),\qquad j=1,\ldots,k. 
\end{align*}
The true value of $\tau =\left(\tau _1,\ldots,\tau _k\right)$ is $\tau_0
=\left(\tau _1\left(\vartheta _0\right),\ldots,\tau _k\left(\vartheta
_0\right)\right)$.  Let us denote the vector MLE $\hat\tau_n=\left(\hat\tau
_{1,n},\ldots,\hat\tau _{k,n}\right)$ and introduce the corresponding Fisher
information matrix $\II_\tau \left(\vartheta _0\right)=\left(\II_\tau
\left(\vartheta _0\right)_{j,i}\right)_{j,i=1,\ldots,k}$
\begin{align*}
\II_\tau \left(\vartheta _0\right)_{j,i}=\int_{ \tau _j  }^{T}\frac{\lambda 
  '\left(t-\tau _j\right)^2}{ \lambda \left(t-\tau
  _j\right)+\lambda _0}{\rm d}t \;\delta _{i,j}
\end{align*} 
where $\delta _{j,i}=\1_{\left\{j=i\right\}}$, i.e., this matrix is diagonal.   
We have the following result
\begin{theorem}
\label{T2} Let the conditions ${\scr R}$ be fulfilled, then the MLE $\hat\tau
_n$ is consistent, asymptotically normal
\begin{align*}
\sqrt{n}\left(\hat\tau_n -\tau _0\right)\Longrightarrow \xi \sim {\cal
  N}\left(0,\II_\tau \left(\vartheta _0\right)^{-1}\right), 
\end{align*}
we have the convergence of polynomial moments  and this MLE is asymptotically
efficient.  
\end{theorem}
{\bf Proof.} The proof of this theorem can be obtained by a slight
modification of the proof of the Theorem 2.4 in \cite{Kut98}. The similar
properties has the corresponding BE $\tilde\tau _n$. 

Suppose that we have already this vector of MLE $\hat \tau _n$ with the values
close to the true value. Our goal is to construct an estimator of the position
$\vartheta _0=\left(x_0,y_0\right)$. We propose the linear system which gives
a good estimator of $\vartheta _0$ as follows. We have
\begin{align*}
z_{j,n}&\equiv \nu ^2\hat\tau _{j,n}^2=\left(x_j-\hat x_0\right)^2+\left(y_j-\hat
y_0\right)^2=x_j^2 +y_j^2+\hat x_0^2+\hat y_0^2-2x_j\hat x_0-2y_j\hat y_0\\
&= r_j^2+\hat r_0^2-2x_j\hat x_0-2y_j\hat y_0=r_j^2-2x_j\hat \gamma
_1-2y_j\hat \gamma _2 +\hat \gamma _3,
\end{align*}
where we denoted $\hat
\gamma _1=\hat x_0,\hat
\gamma _2=\hat y_0,\hat \gamma _3=\|\hat \vartheta _0\|^2
$. Consider the problem of estimation of the vector 
$\gamma_0 =\left(\gamma _{0,1},\gamma _{0,2},\gamma _{0,3}\right)$ by the
observations
\begin{align*}
z_{j,n}=r_j^2-2x_j\gamma _{0,1}-2y_j\gamma _{0,2}+\gamma _{0,3}+\varepsilon
_{j,n},\qquad  j=1,\ldots,k.
\end{align*}
Here $\varepsilon
_{j,n} $ is the noise process and the proposed model is an approximation of
the model of observations. 
Using the method of least squares 
\begin{align*}
\frac{\partial }{\partial \gamma _l
}\sum_{j=1}^{k}\left[z_{j,n}-r_j^2+2x_j\gamma _{1}+2y_j\gamma _{2} -\gamma  
  _{3}\right]^2 =0,\qquad l=1,2,3,
\end{align*}
we obtain the system of equations 
\begin{align*}
&-2\sum_{j=1}^{k}x_j\;\gamma _{1,n}^*-2\sum_{j=1}^{k}y_j\;\gamma _{2,n}^*+k\gamma_{3,n}^*
=\sum_{j=1}^{k}\left(z_{j,n}-r_j^2\right), \\
&-2\sum_{j=1}^{k}x_j^2\;\gamma^*
_{1,n}-2\sum_{j=1}^{k}x_jy_j\;\gamma _{2,n}^* +\sum_{j=1}^{k}x_j\;\gamma_{3,n}^*
=\sum_{j=1}^{k}x_j\left(z_{j,n}-r_j^2\right), \\
&-2\sum_{j=1}^{k}y_jx_j\;\gamma^*
_{1,n}-2\sum_{j=1}^{k}y_j^2\;\gamma _{2,n}^*+\sum_{j=1}^{k}y_j\;\gamma_{3,n}^* 
=\sum_{j=1}^{k}y_j\left(z_{j,n}-r_j^2\right).
\end{align*} 
Therefore we have to solve the linear equation
\begin{align*}
\AA \gamma_n^* =Z_n,
\end{align*}
where
\[\AA=\left(
\begin{array}{ccc}
 -2\sum_{j=1}^{k}x_j,& -2\sum_{j=1}^{k}y_j, &k,\\
-2\sum_{j=1}^{k}x_j^2,&-2\sum_{j=1}^{k}x_jy_j,&\sum_{j=1}^{k}x_j,  \\  
-2\sum_{j=1}^{k}x_jy_j,&-2\sum_{j=1}^{k}y_j^2,&\sum_{j=1}^{k}y_j,  \\  
\end{array}
\right)\]
 and
$Z_n=\left(\sum_{j=1}^{k}\left(z_{j,n}-r_j^2\right),\sum_{j=1}^{k}
x_j\left(z_{j,n}-r_j^2\right),\sum_{j=1}^{k}y_j\left(z_{j,n}-r_j^2\right)
\right)^\TT$.

We consider now the three dimensional parameter $\gamma =\left(\gamma _1,\gamma
_2,\gamma _3\right)^\TT$ with  independent components. Of
course we have the relation $\gamma _3^2=\gamma _1^2+\gamma _2^2$, which we do
not use here. We suppose that the symmetric matrix $\AA$ is non degenerate.

Let us introduce the $k\times 3$ matrix    $\CC=\left(c_{j,r}\right)$  
\begin{align*}
c_{j,1}=2\nu ^2\tau _{0,j}\sigma _j,\quad c_{j,2}=2\nu ^2x_j\tau _{0,j}\sigma
_j,\quad c_{j,3}=2\nu ^2y_j\tau _{0,j}\sigma _j, \quad j=1,\ldots,k,
\end{align*}
 where
\begin{align*}
\sigma _j^2=\left(\int_{\tau _j\left(\vartheta
  _0\right)}^{T}\frac{\lambda '_j\left(t-\tau _j\left(\vartheta
  _0\right)\right)^2}{\lambda_j\left(t-\tau _j\left(\vartheta
  _0\right)\right)+\lambda _0}\; {\rm d}t\right)^{-1},
\end{align*}
and put
\begin{align*}
\DD\left(\vartheta _0\right)=\AA^{-1}\CC^\TT\CC\AA^{-1}.
\end{align*}
The properties of the estimator $\gamma_n ^*=\AA^{-1}Z_n$ are given in the
following proposition.

\begin{proposition}
\label{P1} Let the conditions ${\scr R}$ be fulfilled and the matrix $\AA$ be
non degenerate, then the estimator
$\gamma_n ^* $ is consistent, asymptotically normal
\begin{align*}
\sqrt{n}\left( \gamma_n ^*-\gamma _0\right)\Longrightarrow {\cal N}\left(0,
\DD  \left(\vartheta _0\right)\right),
\end{align*}
and the moments converge. 
\end{proposition}
{\bf Proof.} As we have the asymptotic normality of $\hat\tau_{j,n} $ (Theorem
\ref{T2}) we can
write $\hat\tau_{j,n}=\tau _{0,j}+n^{-1/2}\sigma _j\xi _{j,n}$
and $\xi _{j,n}\Longrightarrow \xi _j$. Recall that $\left(\xi _1,\ldots,\xi
_k\right)$ are i.i.d. ${\cal N}\left(0,1\right)$.
Therefore
\begin{align*}
Z_{1,n}&=\sum_{j=1}^{k}\left(z_{j,n}-r_j^2\right)=\sum_{j=1}^{k}\left(\nu ^2
\hat\tau ^2_{j,n}-r_j^2\right)\\
&=\sum_{j=1}^{k}\left(\nu ^2\left(\tau
_{0,j}+n^{-1/2}\sigma _j\,\xi _{j,n}\right)^2-r_j^2\right)\\
&=
\sum_{j=1}^{k}\left(\nu ^2\tau _{0,j}^2-r_j^2\right) +\frac{2\nu ^2}{\sqrt{n}}
\sum_{j=1}^{k}\tau _{0,j} \sigma _j\,\xi _{j,n}+O\left(\frac{1}{n}\right)\\
&=Z_{1,0}+\frac{1}{\sqrt{n}} \sum_{j=1}^{k} c_{j,1}\,\xi
_{j,n}+O\left(\frac{1}{n}\right) ,
\end{align*}
and
\begin{align*}
Z_{2,n}&=\sum_{j=1}^{k}x_j\left(\nu ^2\tau _{0,j}^2-r_j^2\right) +\frac{2\nu
  ^2}{\sqrt{n}} \sum_{j=1}^{k}x_j\tau _{0,j} \sigma _j\,\xi 
_{j,n}+O\left(\frac{1}{n}\right)\\
&=Z_{2,0}+\frac{1}{\sqrt{n}} \sum_{j=1}^{k} c_{j,2}\,\xi
_{j,n}+O\left(\frac{1}{n}\right), \\
Z_{3,n}&=\sum_{j=1}^{k}y_j\left(\nu ^2\tau _{0,j}^2-r_j^2\right) +\frac{2\nu
  ^2}{\sqrt{n}} 
\sum_{j=1}^{k}y_j\tau _{0,j} \sigma _j\,\xi _{j,n}+O\left(\frac{1}{n}\right)\\
&=Z_{3,0}+\frac{1}{\sqrt{n}} \sum_{j=1}^{k} c_{j,3}\,\xi
_{j,n}+O\left(\frac{1}{n}\right).
\end{align*}
The limit covariance matrix
\begin{align*}
\RR_{l,m}=\lim_{n\rightarrow \infty }n\Ex_{\vartheta _0}
\left(Z_{l,n}-Z_{l,0}\right)\left(Z_{m,n}-Z_{m,0}\right)=\sum_{j=1}^{k}
c_{j,l} c_{j,m} .
\end{align*}
We can write
\begin{align*}
\sqrt{n}\left(\gamma _n^*-\gamma
_0\right)=\AA^{-1}\sqrt{n}\left(Z_n-Z_0\right)\Longrightarrow  \AA^{-1}\CC^\TT \xi  
\end{align*}

The convergence of the moments follow from the convergence of the moments of
the estimator $\hat\tau _n$. 

Let us denote $\vartheta _n^*=\left(\gamma ^*_{1,n},\gamma
^*_{2,n}\right)^\TT$ and call it {\it mean square estimator} (MSE) of $\vartheta $.

Then the Proposition \ref{P1} allows us to write the following
\begin{corollary} Let the conditions ${\scr R}$ be fulfilled and the matrix
  $\AA $ be non degenerate, then the
estimator $\vartheta _n^* $ is consistent and  asymptotically normal
\begin{align*}
\sqrt{n}\left( \vartheta _n^*-\vartheta _0\right)\Longrightarrow {\cal
  N}\left(0, \MM\left(\vartheta _0\right)\right),
\end{align*}
where
\[\MM\left(\vartheta _0\right)=\left(
\begin{array}{cc}
\DD\left(\vartheta _0\right)_{1,1},&\DD\left(\vartheta _0\right)_{1,2}\\
\DD\left(\vartheta _0\right)_{2,1},&\DD\left(\vartheta _0\right)_{2,2} \\  
\end{array}
\right).\]
\end{corollary}
To avoid the large errors we can introduce the following condition
$$
S_n=\left|\gamma _{3,n}^2-\gamma _{1,n}^2-\gamma _{2,n}^2\right|<n^{-1/4}.
$$
For the large values of $n$ it has to be fulfilled. 

\section{The case $\kappa =\frac{1}{2}$}

Let us consider the case of intensity function \eqref{int} in the case $\kappa
=\frac{1}{2}$, i.e.,
\begin{align*}
\lambda _{j,n}\left(\vartheta ,t\right)=an\left|t-\tau _j\right|
^{1/2}\1_{\left\{t\geq \tau _j\right\}}+n\lambda _0,\qquad 0\leq t\leq T ,
\end{align*} 
where $\tau _j=\tau _j\left(\vartheta \right)$ is of course smooth function of
$\vartheta $. Recall that if $\kappa \in
\left(0,\frac{1}{2} \right)$, then we have {\it cusp} case \cite{CDK18} and if
$\kappa >\frac{1}{2}$, then we have the {\it smooth} case considered in this
work. We start with the problem of estimation the parameter $\tau =\left(\tau
_1, \ldots,\tau _k\right)$. Moreover as the Poisson processes with such
intensity functions are independent it is sufficient to study the estimation
of  just one  $\tau _j$, which we denote as $\tau $. Hence we suppose that the
intensity function of the observed Poisson process
$X^n=\left(X\left(t\right)\right),0\leq t\leq T$ is  $\lambda _{n}\left(t-\tau\right)=n\lambda \left(t-\tau\right) $
\begin{align*}
\lambda \left(t-\tau\right)=a\left|t-\tau\right|
^{1/2}\1_{\left\{t\geq \tau \right\}}+\lambda _0,\qquad 0\leq t\leq T .
\end{align*}

Note that  the integral (Fisher information)
\begin{align*}
\II_\tau &=\int_{0}^{T}{ \left(\frac{\partial\lambda \left(t-\tau\right)}{\partial \tau }\right)^2}{\lambda
  \left(t-\tau\right)}^{-1}\,{\rm d}t\\
& =\frac{a^2}{4}\int_{\tau }^{T}\frac{{\rm
    d}t}{
\left|t-\tau \right| \left[ a\left|t-\tau\right|
^{1/2}+\lambda _0 \right] }=\infty .
\end{align*}
Introduce the normalizing function $\varphi _n=\left(n\ln n\right)^{-1/2}$ and
the corresponding log-likelihood ratio process (below $\tau _0$ is the true value
and $u>0$)
\begin{align*}
&\ln Z_n\left(u\right)=\int_{\tau_0}^{T}\ln \frac{\lambda
  \left(t-\tau_0-\varphi _nu\right)}{\lambda
  \left(t-\tau_0\right)} {\rm d}\pi _n\left(t\right)\\
&\qquad -n\int_{\tau_0}^{T}\left[\frac{\lambda
 \left(t-\tau_0-\varphi _nu\right)}{\lambda
  \left(t-\tau_0\right)} -1-   \ln \frac{\lambda
  \left(t-\tau_0-\varphi _nu\right)}{\lambda
  \left(t-\tau_0\right)} \right]\lambda
  \left(t-\tau_0\right){\rm d}t\\
&\qquad \qquad =I_n\left(u\right)-J_n\left(u\right)
\end{align*}
with obvious notations.

We have the asymptotics
\begin{align*}
&\Ex_{\tau _0}I_n\left(u\right)^2=n\int_{\tau_0}^{T}\left(\ln \frac{\lambda
  \left(t-\tau_0-\varphi _nu\right)}{\lambda
  \left(t-\tau_0\right)}\right)^2 \lambda
  \left(t-\tau_0\right) {\rm d}t\\
&=n\int_{\tau_0}^{\tau_0+\varphi _nu}\left( \ln \frac{\lambda _0
  }{a\left(t-\tau _0\right)^{1/2}+ \lambda _0}     \right)^2\lambda
  \left(t-\tau_0\right) {\rm d}t\\
&=n\int_{\tau_0+\varphi _nu}^{T}\left( \ln \frac{a\left(t-\tau _0-\varphi _nu\right)^{1/2}+\lambda _0
  }{a\left(t-\tau _0\right)^{1/2}+ \lambda _0}     \right)^2\lambda
  \left(t-\tau_0\right) {\rm d}t\\
&=o\left(\frac{u}{\ln n}\right)+n\int_{\varphi _nu}^{T-\tau _0}\left( \ln \left(1+\frac{a\sqrt{t-\varphi
      _nu}-a\sqrt{t}  }{a\sqrt{t}+ \lambda _0}  \right)
  \right)^2\left(a\sqrt{t}+ \lambda _0  \right) {\rm d}t\\
&=o\left(\frac{u}{\ln n}\right)+na^2\int_{\varphi _nu}^{T-\tau
    _0}\frac{\left(\sqrt{t-\varphi 
      _nu}-\sqrt{t} \right)^2 }{a\sqrt{t}+ \lambda _0}  
   {\rm d}t\left(1+o\left(1\right)\right).
\end{align*}
Below we  put $t=s\varphi _n$ and use the expansion (for large $s$) 
\begin{align*}
\left(\sqrt{s-u}-\sqrt{s}
\right)^2=s\left(\sqrt{1-\frac{u}{s}}-1\right)^2=\frac{u^2}{4s} \left(1+O\left(\frac{1}{s}\right)\right).
\end{align*}
Therefore
\begin{align*}
&n\int_{\varphi _nu}^{T-\tau _0}\frac{\left(\sqrt{t-\varphi
      _nu}-\sqrt{t} \right)^2 }{a\sqrt{t}+ \lambda _0}  
   {\rm d}t=n\varphi _n^2\int_{u}^{\frac{T-\tau _0}{\varphi
       _n}}\frac{\left(\sqrt{s-u}-\sqrt{s} \right)^2 }{a\sqrt{{s}{\varphi _n}}+
     \lambda _0}      {\rm d}s\\
&\qquad \quad =\frac{n\varphi _n^2u^2}{4}\int_{u}^{\frac{T-\tau _0}{\varphi
       _n}}\frac{1 }{s\left(a\sqrt{{s}{\varphi _n}}+
     \lambda _0\right)}      {\rm d}s\left(1+o\left(1\right)\right)\\
&\qquad \quad =\frac{n\varphi _n^2u^2}{4}\int_{u}^{\frac{T-\tau _0}{\varphi
       _n}}\frac{1 }{s\left(a\sqrt{{s}{\varphi _n}}+
     \lambda _0\right)}      {\rm d}s\left(1+o\left(1\right)\right)\\
&\qquad \quad \approx \frac{n\varphi _n^2u^2}{ 4\left(a\sqrt{T-\tau _0}+
     \lambda _0\right)} \ln \left(\frac{T-\tau _0}{\varphi
       _n} \right)\\
&\qquad \quad \approx \frac{n\varphi _n^2u^2}{4\left(a\sqrt{T-\tau _0}+
     \lambda _0\right) }  \left[\ln \left(n\ln n\right)^{1/2}+\ln  \left(T-\tau
       _0\right)\right]\\
&\qquad \quad\approx\frac{u^2}{8\left(a\sqrt{T-\tau _0}+
     \lambda _0\right)}.
\end{align*}
The stochastic integral admits the representation
\begin{align*}
I_n\left(u\right)&=a\varphi _n\sqrt{n}\int_{u}^{\frac{T-\tau _0}{\varphi
    _n}}\frac{\sqrt{s-u}-\sqrt{s}}{{  a\sqrt{s\varphi _n} +\lambda _0}}\;
{\rm d}W_n\left(s\right) \left(1+o\left(1\right)\right)\\
&=\frac{a\varphi _n\sqrt{n}u}{2}\int_{u}^{\frac{T-\tau _0}{\varphi
    _n}}\frac{1}{\sqrt{s}{\left(  a\sqrt{s\varphi _n} +\lambda _0\right)}}\;
{\rm d}W_n\left(s\right) \left(1+o\left(1\right)\right).
\end{align*}
Here we denoted
\begin{align*}
W_n\left(s\right)&= \frac{1}{\sqrt{n\varphi _n}} \int_{\tau _0}^{\tau
  _0+\varphi _ns}  \left[{\rm 
  d}X\left(v\right)-n \lambda \left(v-\tau _0\right){\rm d}v \right]\\
&=\frac{X\left(\tau _0+s\varphi _n\right)-X\left(\tau
  _0\right)-n\int_{\tau _0}^{\tau _0+s\varphi _n}\lambda \left(t-\tau
  _0\right){\rm d}t }{\sqrt{n\varphi _n}}.
\end{align*}

Using the characteristic function of the stochastic integral we can verify
that  as $n\rightarrow \infty $
\begin{align*}
\frac{1}{\sqrt{\ln n}}\int _u^{\left(T-\tau _0\right)\sqrt{n\ln n}  }\frac{{\rm
    d}W_n\left(s\right)}{\sqrt{s}\left(  a\sqrt{s\varphi _n} +\lambda
  _0\right)}\Longrightarrow  {\cal N}\left(0,\frac{1}{ a\sqrt{T-\tau
    _0}+\lambda _0}\right). 
\end{align*}
The similar calculations for the ordinary integral $J_n\left(u\right)$ give
us the asymptotics
\begin{align*}
J_n\left(u\right)=\frac{u^2a^2}{16\left(a\sqrt{T-\tau _0}+
     \lambda _0\right)}+o\left(1\right).
\end{align*}
Theefore for the likelihood ratio process  we have
\begin{align}
\label{qq}
Z_n\left(u\right)=\exp\left\{\gamma u\Delta _n -\frac{\gamma ^2u^2}{2}
+r_n\right\} 
\end{align}
where $ r_n\rightarrow 0$,
\begin{align}
\label{qk}
\gamma ^2=\frac{a^2}{8\left(a\sqrt{T-\tau _0}+\lambda _0\right)},\qquad \Delta 
_n\Longrightarrow  {\cal N}\left(0,1\right) .
\end{align}
The  estimates \eqref{sm} and \eqref{pi} can be obtained too. These estimates 
and the representation \eqref{qq}-\eqref{qk} allow verify the convergences
\begin{align*}
\sqrt{n\ln n}\left(\hat\tau _n-\tau _0\right)\Longrightarrow {\cal
  N}\left(0,\gamma ^{-2}\right) 
\end{align*}
and
\begin{align*}
 \Ex_{\tau _0}\left(\hat\tau _n-\tau _0\right)^2=\frac{1}{\gamma ^{2}n\ln
   n}\left(1+o\left(1\right)\right). 
\end{align*}
Therefore this is {\it smooth} or regular case with asymptotically normal MLE.

\section{One-step MLE}

We have the same model of observations but our problem is to construct an
estimator-process $\vartheta ^\star=(\vartheta ^\star_{t,n}, 0\leq t\leq T)$,
where the estimator $\vartheta ^\star_{t,n}$ has such form that it can be easy
calculated and is asymptotically efficient. Of course we cannot use the MLE
because its calculation for all $t\in(0,T]$ is computationally too difficult.

Such construction in the case of ergodic diffusion processes was proposed in
the work \cite{K17}. See as well the work \cite{KhK18}, where the similar
approach was applied in the case of parameter estimation of hidden telegraph
process.  The case of inhomogeneous Poisson processes was
considered in \cite{DGK18}. Here we apply the developed there techniques. The
main advantage of this approach is the simplicity of calculations of
asymptotically efficient estimators.

We need  a consistent and asymptotically normal estimator 
$\vartheta _n^*$ 
\begin{align}
\label{mse}
n^{q}\left(\vartheta _n^*-\vartheta _0\right)\Longrightarrow {\cal N}\left(0,
\MM\left(\vartheta _0\right)\right) ,
\end{align}
where $q\in \left(\frac{1}{4},\frac{1}{2}\right)$ and $\MM\left(\vartheta
_0\right) $ is some non degenerate matrix.  To construct such estimator we
follow the work \cite{Kh09}, Section 3.3 and use the thinning of a Poisson
process. Let $X^T=\left(X\left(t\right),0\leq t\leq T\right)$ be the Poisson
process with intensity function $\lambda \left(t\right)$. Let
$t_i,i=1,2,\ldots $ be the events of the process $X^T$, so that
\begin{align*}
X\left(t\right)=\sum_{i}^{}\1_{\left\{t_i<t\right\}}.
\end{align*}
Let $\eta _1,\eta _2,\ldots$ be i.i.d.  random variables such
that $\Pb\left(\eta _i=1\right)=p=1-\Pb\left(\eta _i=0\right)$. Introduce the
new process $Y\left(t\right)$ as follows
\begin{align*}
Y\left(t\right)=\sum_{i}^{}\eta _i\1_{\left\{t_i<t\right\}},\qquad 0\leq t\leq T.
\end{align*}
\begin{lemma} (\cite{R06}, Proposition 5.2)
\label{L4}
The process $Y\left(t\right),0\leq t\leq T$ is a Poisson process with the
intensity function $p\lambda \left(t\right)$. The process $\tilde
X\left(t\right)=X\left(t\right)-Y\left(t\right) $ is also Poisson with
intensity function $\left(1-p\right)\lambda \left(t\right)$. The processes
$Y\left(t\right)$ and $\tilde X\left(t\right)$ are independent.

\end{lemma} Note that in \cite{R06} the proof is given  for Poisson processes
with constant intensity function, but it can be easily modified to cover
inhomogeneous Poisson processes considered in our work. 

The observed $k$ Poisson processes $X^n=\left(X_j\left(t\right),0\leq t\leq T,
j=1,\ldots,k\right)$   with intensity functions
\begin{align*}
\lambda _{j,n}\left(\vartheta _0\right)=n\lambda_j \left(t-\tau_j
\left(\vartheta _0\right)\right) +n\lambda _0,\quad 0\leq t\leq T,\quad j=1,\ldots,k
\end{align*}
using the thinning procedure we represent as the sum of $2k$ independent
Poisson processes 
 \begin{align*}
Y^n& =\left(Y_j\left(t\right),0\leq t\leq T, j=1,\ldots,k\right) ,\\
\tilde X^n& =(\tilde
X_j\left(t\right),0\leq t\leq T, j=1,\ldots,k),
\end{align*}
where $Y_j\left(t\right)=
X_j\left(t\right)-\tilde
X_j\left(t\right)$ with the intensity functions
\begin{align*}
\lambda _{j,n}^{Y}\left(\vartheta _0,t\right)&=np_n\lambda_j \left(t-\tau_j
\left(\vartheta _0\right)\right) +np_n\lambda _0,\quad\qquad  0\leq t\leq T,\\ 
\lambda _{j,n}^{\tilde X}\left(\vartheta
_0,t\right)&=n\left(1-p_n\right)\lambda_j \left(t-\tau_j \left(\vartheta
_0\right)\right) +n\left(1-p_n\right)\lambda _0,\quad 0\leq t\leq T,
\end{align*}
respectively. We put the probability $p_n=n^{-b}$, where $b\in
\left(0,\frac{1}{2}\right)$. 

Let us denote $\vartheta _n^*$ the MSE constructed by the observations $Y^n$
and remark that it is asymptotically normal \eqref{mse} with
$q=\frac{1-b}{2}$. Remark that we need not to use all $k$ detectors and it is
sufficient to construct the preliminary estimator on the base of the three
Poisson processes from three detectors (not on the same line). 
 We can not introduce the One-step MSE-process $\vartheta _{t,n}^\star, t\in
 (0,T] $   as in \cite{DGK18} by the formula 
\begin{align*}
\vartheta _{t,n}^\star&=\vartheta _n^* +\II_t\left(\vartheta _n^*\right)^{-1}
\sum_{j=1}^{k}\frac{\partial \tau_j\left(\vartheta _n^*\right) }{\partial
  \vartheta } \int_{\tau_j\left(\vartheta _n^*\right)}^{t}\frac{\ell_j\left(s,\vartheta
  _n^*\right)}{{n}}\left[{\rm 
    d}\tilde X_j\left(s\right)-\lambda _{j,n}^{\tilde X}\left(\vartheta
  _n^*,s\right){\rm d}s\right].
\end{align*}
Here we suppose that $\frac{0}{0}=0$, the function
\begin{align*}
\ell_j\left(s,\vartheta _n^*\right)=\frac{\lambda_j '\left(s-\tau
  _j\left(\vartheta _n^* \right)\right)}{\lambda_j\left(s-\tau
  _j\left(\vartheta _n^* \right)\right)+\lambda _0}\,\1_{\left\{s>\tau
  _j\left(\vartheta _n^*\right)\right\}} 
\end{align*}
and the Fisher information matrix $\II_t\left(\vartheta\right) $ is the
slightly modified matrix
$\II\left(\vartheta\right) $. The modification concerns the
weights $J _j\left(\vartheta \right)$ in the definition of the norm
$\left\|a\right\|_\vartheta $ and scalar product $\langle a,b\rangle_\vartheta
$. The modified weights are 
\begin{align*}
J_{j,t}\left(\vartheta \right) =\frac{1}{\nu ^2\left\|\vartheta_j
  -\vartheta\right\|}\int_{\tau _j\left(\vartheta 
  \right)}^{t}\frac{\lambda_j '\left(s-\tau _j\left(\vartheta
   \right)\right)^2}{\lambda_j\left(s-\tau _j\left(\vartheta 
  \right)\right)+\lambda _0}\;\1_{\left\{s>\tau _j\left(\vartheta 
  \right)\right\}}\;{\rm d}s
\end{align*}
and we write $\langle a,b\rangle_{t,\vartheta} $, $\left\|a\right\|_{t,\vartheta}
$. Therefore
\[\II_t\left(\vartheta \right)=\left(
\begin{array}{cc}
\left\|x-x _0\right\|_{t,\vartheta} ^2,&\langle
\left(x-x_0\right),\left(y-y_0\right)\rangle_{t,\vartheta}\\
\langle \left(x-x_0\right),\left(y-y_0\right)\rangle_{t,\vartheta},&\left\|y-y
_0\right\|_{t,\vartheta} ^2 \\  
\end{array}
\right).\]

Let us put the estimators $\tau _j\left(\vartheta _n^*\right)$ in the order of
increasing
$$
\tau _{\left(1\right)}\left(\vartheta _n^*\right)<\tau
_{\left(2\right)}\left(\vartheta _n^*\right)<  \ldots <\tau
_{\left(k\right)}\left(\vartheta _n^*\right).
$$ It is evident that on the time interval $\left[0,\tau
  _{\left(1\right)}\left(\vartheta _n^*\right)\right]$ we have $\vartheta
_{t,n}^\star=\vartheta _n^*$. Moreover, for the values $\tau
_{\left(1\right)}\left(\vartheta _n^*\right)\leq t\leq \tau
_{\left(2\right)}\left(\vartheta _n^*\right) $ the Fisher information matrix
$\II_t\left(\vartheta_n^* \right) $ is degenerated. In the case $\tau
_{\left(2\right)}\left(\vartheta _n^*\right)\leq t\leq \tau
_{\left(3\right)}\left(\vartheta _n^*\right) $ this matrix is non degenerate
and the estimator $\vartheta _{t,n}^\star$ is asymptotically normal.

 Note as well that in the stochastic integral used in the definition
of the One-step MLE-process  the random vector $\vartheta _n^*$  is
independent of the ``observations'' $\tilde X_j\left(\cdot \right)$ because
the Poisson processes $Y^n$ and $\tilde X^n$ are independent. Therefore the
stochastic integral is well defined.

We do not give here the strict proofs but just show why the estimator-process
$\vartheta _{t,n}^\star,0<t\leq T $ is asymptotically normal with the same
parameters as the MLE $\hat\vartheta _{t,n}, 0<t\leq T$. Here $\hat\vartheta
_{t,n}$ is the MLE constructed by the first observations
$X^{t,n}=(X_j\left(s\right),0\leq s\leq t,$ $ j=1,\ldots,k)$. Of course, this
MLE is not even consistent for the values $t\in [0,\tau _{\left(3\right)}] $
because up to $\tau _{\left(1\right)} $ the observations $X^{t,n}$ do not
contain any information about $\vartheta _0$. The consistency of it is
possible for the values $t>\tau _{\left(3\right)} $ only. Indeed for these
values we have the consistent estimators of $\tau
_{\left(1\right)}\left(\vartheta _0\right),\tau
_{\left(2\right)}\left(\vartheta _0\right),\tau
_{\left(3\right)}\left(\vartheta _0\right) $ and if the corresponding
detectors are not on the same line then $\vartheta _0$ can be identified.  The
proof of this asymptotic equivalence is a slight modification of the proof
given in \cite{DGK18}.

We can write
\begin{align*}
&\sqrt{n}\left(\vartheta _{t,n}^\star-\vartheta
  _0\right)=\sqrt{n}\left(\vartheta _n^*-\vartheta _0\right)\\ &\quad
  +\II_t\left(\vartheta _n^*\right)^{-1} \sum_{j=1}^{k}\frac{\partial
    \tau_j\left(\vartheta _n^*\right) }{\partial \vartheta }
  \int_{\tau_j\left(\vartheta _n^*\right) }^{t}\frac{\ell_j\left(s,\vartheta
    _n^*\right)}{\sqrt{n}}{\rm d}\tilde\pi _j\left(s\right)\\ &\quad
  +\II_t\left(\vartheta _n^*\right)^{-1} \sum_{j=1}^{k}\frac{\partial
    \tau_j\left(\vartheta _n^*\right) }{\partial \vartheta }
  \int_{\tau_j\left(\vartheta _n^*\right) }^{t}\frac{\ell_j\left(s,\vartheta
    _n^*\right)}{\sqrt{n}}\left[\lambda_{j,n}^{\tilde X}\left(\vartheta _0,s
    \right)-\lambda_{j,n}^{\tilde X}\left(\vartheta _n^*,s \right)\right]{\rm
    d}s\\ &\quad =\sqrt{n}\left(\vartheta _n^*-\vartheta _0\right)
  +\II_t\left(\vartheta _0\right)^{-1} \sum_{j=1}^{k}\frac{\partial
    \tau_j\left(\vartheta _0\right) }{\partial \vartheta }
  \int_{\tau_j\left(\vartheta _n^*\right) }^{t}\frac{\ell_j\left(s,\vartheta
    _0\right)}{\sqrt{n}}{\rm d}\tilde\pi
  _j\left(s\right)+o\left(1\right)\\ &\; -\II_t\left(\vartheta _0\right)^{-1}
  \sum_{j=1}^{k}\frac{\partial \tau_j\left(\vartheta _0\right) }{\partial
    \vartheta }\frac{\partial \tau_j\left(\vartheta _0\right) }{\partial
    \vartheta }^\TT \int_{\tau_j\left(\vartheta _n^*\right)
  }^{t}\frac{\lambda_j'\left(s-\tau _j\left(\vartheta _0
    \right)\right)^2}{\lambda_j\left(s-\tau _j\left(\vartheta _0
    \right)\right)+\lambda _0 }{\rm d}s\sqrt{n}\left(\vartheta _n^*-\vartheta
  _0\right)\\ &\quad =\II_t\left(\vartheta _0\right)^{-1}
  \sum_{j=1}^{k}\frac{\partial \tau_j\left(\vartheta _0\right) }{\partial
    \vartheta }\frac{\partial \tau_j\left(\vartheta _0\right) }{\partial
    \vartheta }^\TT \frac{1}{\sqrt{n}}\int_{\tau_j\left(\vartheta
    _n^*\right)}^{t}{\ell_j\left(s,\vartheta _0\right)}{\rm d}\tilde\pi
  _j\left(s\right)+o\left(1\right).
\end{align*}
Here ${\rm d}\tilde\pi _j\left(s\right)={\rm d}\tilde
X_j\left(s\right)-\lambda_{j,n}^{\tilde X}\left(\vartheta _0,s \right){\rm d}s$, we
used the consistency of the estimator $\vartheta _n^*=\vartheta
_0+O\left(\frac{1}{\sqrt{n^{1-b}}}\right) $, Taylor formula and the equality
\begin{align*}
\II_t\left(\vartheta _n^*\right)=\sum_{j=1}^{k}\frac{\partial
    \tau_j\left(\vartheta _0\right) }{\partial \vartheta }\frac{\partial
    \tau_j\left(\vartheta _0\right) }{\partial \vartheta }^\TT \int_{
  \tau_j\left(\vartheta _0\right) 
  }^{t}\frac{\lambda_j'\left(s-\tau
  _j\left(\vartheta _0 \right)\right)^2}{\lambda_j\left(s-\tau
  _j\left(\vartheta _0 \right)\right)+\lambda _0 }{\rm
  d}s+O\left(\frac{1}{\sqrt{n^{1-b}}} \right).
\end{align*}
Hence (see \eqref{mse})
\begin{align*}
\sqrt{n}\left(\vartheta _n^*-\vartheta _0\right)O\left(n^{-\frac{1-b}{2}}
\right)=n^q \left(\vartheta _n^*-\vartheta
_0\right)O\left(n^{\frac{1}{2}-q-\frac{1-b}{2}}\right)=n^q \left(\vartheta
_n^*-\vartheta _0\right) o\left(1\right).
\end{align*}
 Now the asymptotic normality
\begin{align*}
\sqrt{n}\left(\vartheta _{t,n}^\star-\vartheta _0\right)\Longrightarrow {\cal
  N}\left(0,\II_t\left(\vartheta _0\right)^{-1}\right)
\end{align*}
follows from the central limit theorem for (independent) stochastic integrals 
\begin{align*}
\frac{1}{\sqrt{n}}\int_{0 }^{t}{\ell_j\left(s,\vartheta
    _0\right)}{\rm d}\tilde \pi _j\left(s\right),\qquad j=1,\ldots,k.
\end{align*}
Therefore, the One-step MLE-process $\vartheta _{t,n}^\star$ is asymptotically
equivalent to the MLE $\hat\vartheta _{t,n}$. Of course, the it is possible to
verify the uniform convergence of moments and the asymptotic efficiency of
this estimator.

Recall that the One-step MLE-process is  uniformly consistent, i.e.,
for any $\varepsilon >0$ and any $\nu >0$
\begin{align*}
\Pb_{\vartheta _0}\left(\sup_{\varepsilon \leq t\leq T}\left\| \vartheta
_{t,n}^\star-\vartheta _0 \right\|>\nu  \right)\longrightarrow 0.
\end{align*}
Moreover, it is possible to verify the weak convergence of the random  process
$ u_{t,n}^\circ =\sqrt{n}\left(\vartheta _{t,n}^\star-\vartheta _0 \right),
\varepsilon \leq t\leq T $ to the limit Gaussian process (see \cite{DGK18}).

Of course, we can have the One-step MLE $\vartheta _{n}^\star=\vartheta
_{T,n}^\star$ too
\begin{align*}
\vartheta _{n}^\star&=\vartheta _n^* +\II\left(\vartheta _n^*\right)^{-1}
\sum_{j=1}^{k}\frac{\partial \tau_j\left(\vartheta _n^*\right) }{\partial
  \vartheta } \int_{\tau_j\left(\vartheta _n^*\right) }^{T}\frac{\ell_j\left(t,\vartheta
  _n^*\right)}{{n}}\left[{\rm d}\tilde X_j\left(t\right)-\lambda
  _{j,n}^{\tilde X}\left(\vartheta _n^*,t\right){\rm d}t\right].
\end{align*}
This estimator have the same asymptotic normality
\begin{align*}
\sqrt{n}\left(\vartheta _{n}^\star-\vartheta _0\right)\Longrightarrow {\cal
  N}\left(0,\II_T\left(\vartheta _0\right)^{-1}\right).
\end{align*}
We can check the uniform convergence of the moments (see \cite{DGK18}) and
therefore to prove the asymptotic efficiency of this estimator. 

Another possibility is to use the CUSUM type estimators of $\tau _j$ and then
having first three estimators of the closest to the source sensors we can
start the  One-step MLE of the position of source based on the observations of
other sensors.

\section{Discussion}

In this work we supposed that the source starts  emission  at the
instant $t=0$. It is interesting to consider the more general statement with
the unknown beginning of emission $\tau_*$. Therefore the signal   received 
by the $j$-th detector arrives at the moment $\bar\tau _j=\tau ^*+\tau _j $,
where $\tau _j=\nu ^{-1}\left\|\vartheta _j-\vartheta _0\right\|$. Let us
denote $\hat\tau _{j,n}$ the MLE of the arriving time in the $j$-th
detector. Then we have
\begin{align*}
\nu ^2\left(\bar\tau _j-\tau_*\right)^2=\left(x_j-x_0\right)^2
+\left(y_j-y_0\right)^2 ,\qquad j=1,\ldots,k,
\end{align*}
and
\begin{align*}
\nu ^2\bar\tau _{j}^2=x_j^2+y_j^2+x_0^2+y_0^2- \nu ^2\tau
_*^2-2x_jx_0-2y_jy_0+2\nu ^2\bar\tau _{j} \tau _*.
\end{align*}
Let us denote $\gamma _1=x_0,\gamma _2=y_0,\gamma _3=\tau _*,\gamma
_4=x_0^2+y_0^2- \nu ^2\tau _*^2$ and $r_j^2=x_j^2+y_j^2$. Therefore for the
estimator $\gamma _n^*=\left(\gamma _{1,n}^*,\ldots,\gamma _{4,n}^*\right)$ of
the parameter $\gamma =\left(\gamma _1,\ldots,\gamma _4\right)$ we obtain 
the system
\begin{align*}
-2x_j\gamma _{1,n}^*-2y_j\gamma _{2,n}^*+2\nu ^2\hat\tau _{j,n}\gamma
_{3,n}^*+\gamma _{4,n}^*=\nu ^2\hat\tau _{j,n}^2 -r_j^2,\qquad j=1,\ldots,k
\end{align*}
and so on. The corresponding matrix $\AA$ is already random and the estimator
needs a special study.  From the consistency of the estimators $\hat\tau
_{j,n},j=1,\ldots,k $ we obtain the consistency of the estimator $\gamma
_n^*$.

{\bf Acknowledgment.} This work was done under partial financial support of
the grant of RSF number 14-49-00079.

\end{document}